\documentclass[12pt,twoside]{amsart}
\usepackage{amssymb}
\usepackage{verbatim}
\usepackage{amsmath}
\usepackage{bm}
\usepackage{a4wide}
\usepackage[utf8]{inputenc}
\usepackage[T1]{fontenc}
\usepackage{times}
\usepackage{latexsym}
\usepackage{enumerate}
\usepackage{needspace}
\usepackage[hidelinks]{hyperref}

\newcommand{\R}{\ensuremath{\mathbb{R}}}
\newcommand{\D}{\ensuremath{\mathbb{D}}}
\newcommand{\C}{\ensuremath{\mathbb{C}}}

\makeatletter
\newcommand{\sumprime}{\if@display\sideset{}{'}\sum\else\sum'\fi}
\makeatother
\allowdisplaybreaks[2]
\begin{document}
\numberwithin{equation}{section}
\newtheorem{theorem}{Theorem}[section]
\newtheorem{proposition}[theorem]{Proposition}
\newtheorem{conjecture}[theorem]{Conjecture}
\def\theconjecture{\unskip}
\newtheorem{corollary}[theorem]{Corollary}
\newtheorem{lemma}[theorem]{Lemma}
\newtheorem{observation}[theorem]{Observation}
\newtheorem{definition}{Definition}
\numberwithin{definition}{section}
\theoremstyle{remark}
\newtheorem{remark}[theorem]{Remark}
\newtheorem*{remark*}{Remark}
\theoremstyle{plain}
\newtheorem{kl}{Key Lemma}
\def\thekl{\unskip}
\newtheorem{question}{Question}
\def\thequestion{\unskip}
\newtheorem{example}{Example}
\def\theexample{\unskip}
\newtheorem{problem}{Problem}

\thanks{The author is supported by the National Natural Science Foundation of Sichuan, No. 2025ZNSFSC0796, the National Natural Science Foundation of China, No. 12401103, and the Fundamental Research Fund for the Central Universities, Southwest Minzu University, No. ZYN2025007.}
\title[Bergman--Szeg\H{o} comparison]
{Comparison of the Bergman and Szeg\H{o} kernels on arbitrary bounded
pseudoconvex domains with $C^2$ boundary}
\author[Xu Xing]{Xu Xing}
\begin{abstract}
On every bounded pseudoconvex domain $\Omega\Subset\C^n$, $n\ge2$,
with $C^2$ boundary, we prove
$K_\Omega/S_\Omega\le C_{\Omega,\varepsilon}
\delta_\Omega^{-1}|\log\delta_\Omega|^{1+\varepsilon}$
near the boundary for every $\varepsilon>0$,
where $\delta_\Omega$ denotes the boundary distance.
For each positive nondecreasing $H(s)=o(s)$, a smoothly bounded
pseudoconvex complete Reinhardt domain in $\C^2$ shows that the logarithmic factor
cannot be replaced by $H(|\log\delta_\Omega|)$.
\bigskip


\smallskip

\end{abstract}

\address[Xu Xing]{School of Mathematics, Southwest Minzu University, Chengdu, 610041, China
\newline\mbox{\quad}\,Department of Mathematics and Computer Science,
\newline\mbox{\quad}\,University of Cologne, Cologne, 50931, Germany}
\email{xingxu@swun.edu.cn}
\maketitle

\section{Introduction}
The Bergman and Szeg\H{o} kernels are fundamental objects linking
holomorphic function theory, boundary geometry, and partial differential
equations. Classical work on kernel functions
and conformal mapping revealed their connections with extremal
problems and potential theory; see Bergman--Schiffer \cite{BergmanSchiffer51}, Stein \cite{Stein}, and
Kerzman--Stein \cite{KerzmanStein1978}. In several complex variables,
the Bergman kernel also determines a biholomorphically invariant metric, while the Szeg\H{o}
projection provides a basic analytic tool for studying holomorphic boundary
values and the tangential Cauchy--Riemann equations. Estimates for these
kernels and their projections play important roles in the boundary regularity
of holomorphic mappings, $\bar\partial$-Neumann theory, and CR analysis
\cite{BellLigocka80,BS99,BdMS76}.

For a bounded domain $\Omega\subset\C^n$, let
$\delta_\Omega(z)=\operatorname{dist}(z,\partial\Omega)$.
We denote the Bergman kernel for Euclidean volume by $K_\Omega(z,w)$ and,
when $\partial\Omega$ is $C^2$, the Szeg\H{o} kernel for Euclidean surface
measure by $S_\Omega(z,w)$. Their diagonal values are denoted by
\[
K_\Omega(z):=K_\Omega(z,z),\quad
S_\Omega(z):=S_\Omega(z,z).
\]

On smoothly bounded strongly pseudoconvex domains, the landmark
expansions of Fefferman \cite{Fefferman74} and
Boutet de Monvel--Sj\"ostrand \cite{BdMS76} describe the boundary
singularities of these kernels. The leading terms encode the Levi
geometry, linking complex geometry with microlocal analysis. In particular,
\[
K_\Omega(z)\asymp\delta_\Omega(z)^{-n-1},\quad
S_\Omega(z)\asymp\delta_\Omega(z)^{-n},\quad
\frac{K_\Omega(z)}{S_\Omega(z)}\asymp\delta_\Omega(z)^{-1}.
\]
When the Levi form degenerates, the geometry of complex tangential
contact becomes decisive. Catlin's finite-type subelliptic theory
\cite{Catlin87} and the estimates of Nagel--Rosay--Stein--Wainger in
$\C^2$ \cite{NRSW89} relate analytic bounds to anisotropic geometric
scales. More recent work of Hsiao--Savale \cite{HS22} develops
Bergman--Szeg\H{o} asymptotics on weakly pseudoconvex finite-type domains in $\C^2$.

At infinite type, smoothness of the boundary alone does not ensure
Sobolev regularity of the Bergman projection, as Barrett's analysis of the
Diederich--Forn\ae ss worm demonstrates \cite{Barrett92}.
Nevertheless, sharp kernel estimates can still be obtained under
additional structure; examples include generalized decoupled domains
with controlled infinite-type profiles \cite{Jaiswal25}.
For general weakly pseudoconvex domains, estimates in terms of boundary
distance and plurisubharmonic exhaustion functions offer a complementary
approach \cite{Chen17,Chen23,CZ22}.

The comparison of the Bergman and Szeg\H{o} kernels was proposed by
Stein \cite[p.~20]{Stein}. It asks whether the strongly pseudoconvex
comparison persists when the Levi form degenerates.
On bounded pseudoconvex domains with $C^2$ boundary satisfying
$\delta$-regularity, Chen--Fu \cite{CF} proved the upper estimate
\[
\frac{K_\Omega(z)}{S_\Omega(z)}
\le\frac{C_{\Omega,\eta}|\log\delta_\Omega(z)|^{1/\eta}}
{\delta_\Omega(z)}
\]
for a suitable Diederich--Forn\ae ss exponent $0<\eta\le1$. They also proved the logarithm-free
comparison on bounded convex domains with $C^2$ boundary.
Our first result removes the $\delta$-regularity assumption and
establishes the upper comparison for every logarithmic exponent
greater than one.

\begin{theorem}\label{thm:unweighted}
Let $\Omega\Subset\C^n$, $n\ge2$, be a bounded pseudoconvex domain
with $C^2$ boundary. For every $0<\alpha<1$, there is a constant
$C_{\Omega,\alpha}>0$ such that
\begin{equation}\label{intro:ratio-parameters}
 \frac{K_\Omega(z)}{S_\Omega(z)}
 \le\frac{C_{\Omega,\alpha}|\log\delta_\Omega(z)|^{1/\alpha}}{\delta_\Omega(z)}
\end{equation}
near the boundary.
\end{theorem}
In one complex dimension, the classical theory of planar
kernel functions yields
\[
\lim_{z\to\partial\Omega}
\delta_\Omega(z)\frac{K_\Omega(z)}{S_\Omega(z)}=1
\]
for every bounded domain $\Omega\subset\C$ with $C^2$ boundary.
For related results on planar kernel functions and their identities, see
\cite{KerzmanStein1978,BellGustafsson2022}.
The comparison also yields information about the Szeg\H{o} kernel itself.
The Ohsawa--Takegoshi extension theorem and Fu's estimate
\cite{OT,Fu94} give $K_\Omega(z)\ge c_\Omega\delta_\Omega(z)^{-2}$
on smoothly bounded pseudoconvex domains. Chen extended this lower
bound to pseudoconvex domains with $C^0$ boundary
\cite{Chen23}.
Theorem~\ref{thm:unweighted} therefore implies
\[
S_\Omega(z)\ge
\frac{c_{\Omega,\alpha}}
{\delta_\Omega(z)|\log\delta_\Omega(z)|^{1/\alpha}},
\quad 0<\alpha<1.
\]
Using Shaw's $L^2$ existence theorem for $\bar\partial_b$,
\cite[Theorem~1.3]{WX} establishes
$S_\Omega(z)\ge c|z-p|^{-1/3}$ for non-tangential approach to every
$p\in\partial\Omega$ on smoothly bounded pseudoconvex domains
$\Omega\subset\C^2$.
The logarithm-free lower bound remains a separate question
\cite[Problem~2]{WX}.

\refstepcounter{theorem}
\begingroup
\renewcommand{\thequestion}{\thetheorem}
\begin{question}[Ohsawa]\label{qu:ohsawa}
Does every bounded pseudoconvex domain $\Omega\Subset\C^n$, $n\ge2$,
with $C^2$ boundary satisfy
\[
 S_\Omega(z)\ge c_\Omega\delta_\Omega(z)^{-1}
\]
near its boundary?
\end{question}
\endgroup

Distance weights provide a natural interpolation between interior and
boundary norms. Let $K_{\Omega,t}(z,w)$ be the reproducing kernel for
$\delta_\Omega^t\,dV$, $t>-1$, and write
\[
K_{\Omega,t}(z):=K_{\Omega,t}(z,z),\quad
K_{\Omega,0}(z)=K_\Omega(z).
\]
Weighted kernels and holomorphic Sobolev scales have also been studied
in the strongly pseudoconvex setting \cite{Englis08}.
Here the distance-weighted $\bar\partial$ estimate in
\cite[Theorem~$1.1'$]{Chen14} allows the comparison to extend to every
$t>-1$, with constants locally uniform in the weight parameter.

\begin{theorem}\label{thm:comparison}
Let $\Omega\Subset\C^n$, $n\ge2$, be a bounded pseudoconvex domain
with $C^2$ boundary. For every $0<\alpha<1$ and
$I\Subset(-1,\infty)$, there is a constant $C_I>0$ such that,
uniformly for $t\in I$,
\begin{equation}\label{intro:weighted}
 \frac{K_{\Omega,t}(z)}{S_\Omega(z)}
 \le C_I\left[\frac{|\log\delta_\Omega(z)|^{1/\alpha}}{\delta_\Omega(z)}\right]^{1+t}
\end{equation}
near the boundary. The constant may depend on $\Omega$, $\alpha$,
and $I$.
\end{theorem}

The upper logarithmic loss is nearly optimal in the following sense.
Proposition~\ref{ex:obstructions-1} constructs, for each prescribed
positive nondecreasing function $H(s)=o(s)$, a smoothly bounded
pseudoconvex complete Reinhardt domain in $\C^2$ on which
\[
\frac{K_\Omega(z)}{S_\Omega(z)}
\le\frac{C_\Omega H(|\log\delta_\Omega(z)|)}{\delta_\Omega(z)}
\]
cannot hold near the boundary. In particular, the general estimate
$K_\Omega/S_\Omega\lesssim\delta_\Omega^{-1}$ is false. Consequently, a logarithm-free upper comparison combined with
$K_\Omega\gtrsim\delta_\Omega^{-2}$ cannot provide a general solution to
Question~\ref{qu:ohsawa}.

The proof combines Chen's distance-weighted $L^2$ estimate
\cite{Chen14}, pluricomplex Green functions \cite{BlockiNotes},
and Diederich--Forn\ae ss barriers \cite{DF}.
The key to removing $\delta$-regularity from the upper comparison is
the weighted kernel interpolation inequality
\[
K_{\Omega,t}(z)
\le C_{\Omega,t}K_\Omega(z)^{1+t}S_\Omega(z)^{-t},
\quad -1<t<0.
\]
Proposition~\ref{prop:kernel-inequality} proves this estimate directly
from a holomorphic decomposition into a Bergman part and a Hardy
part. The two norms are controlled by a weighted Bergman projection
and the multiplicative Sobolev trace inequality; balancing their
pointwise estimates gives the comparison above. The required Sobolev
regularity is supplied for $C^2$ boundaries by Harrington
\cite{Harrington09}, following Kohn's smooth-boundary method
\cite{Kohn73}.
At $t=-1/2$, the local barriers and weighted extension estimates give
\[
c_{\Omega,\alpha}
\frac{\delta_\Omega(z)^{1/2}}
{|\log\delta_\Omega(z)|^{1/(2\alpha)}}K_\Omega(z)
\le K_{\Omega,-1/2}(z)
\le C_\Omega\bigl(K_\Omega(z)S_\Omega(z)\bigr)^{1/2}.
\]
Combining the two inequalities yields the upper comparison under
pseudoconvexity and $C^2$ boundary regularity alone. The same weighted
extension estimates then give Theorem~\ref{thm:comparison} for all
$t>-1$.

Section~\ref{sec:prelim} recalls the analytic tools and Green functions.
Section~\ref{sec:two-sided} proves the weighted upper comparison
by means of the direct interpolation estimate in
Proposition~\ref{prop:kernel-inequality}.
Section~\ref{sec:counterexamples} constructs the logarithmic
obstruction on a Hartogs domain.

\section{Preliminaries}\label{sec:prelim}
We write $B(a,r)=\{z:|z-a|<r\}$ for a Euclidean ball.
The measures $dV$ and $d\sigma$ are Euclidean volume and surface
measure, respectively. As usual,
$X\asymp Y$ means $cY\le X\le CY$ for positive constants independent
of the variables under consideration. Unless stated otherwise, constants
may change from one occurrence to the next and depend on the domain
and fixed parameters.

\subsection{Holomorphic function spaces}\label{sec:holomorphic-prelim}
Let $\Omega\subset\C^n$ be a bounded domain, and let
$\varphi\in L^\infty_{\mathrm{loc}}(\Omega)$ be real-valued.
We define the weighted Bergman space by
\[
 A_\varphi^2(\Omega)
 =\left\{f\in\mathcal O(\Omega):
 \|f\|_\varphi^2:=\int_\Omega |f|^2e^{-\varphi}\,dV<\infty\right\}.
\]
This is a closed subspace of $L^2(\Omega,e^{-\varphi}dV)$.
Its reproducing kernel $K_{\Omega,\varphi}(z,w)$ satisfies
\[
 f(z)=\int_\Omega K_{\Omega,\varphi}(z,w)
 f(w)e^{-\varphi(w)}\,dV(w),
 \quad f\in A_\varphi^2(\Omega).
\]
Writing $K_{\Omega,\varphi}(z):=K_{\Omega,\varphi}(z,z)$, we have
\[
 K_{\Omega,\varphi}(z)
 =\sup\left\{|f(z)|^2:
 f\in A_\varphi^2(\Omega),\ \|f\|_\varphi\le1\right\}.
\]
The orthogonal Bergman projection is
\[
 P_{\Omega,\varphi}:
 L^2(\Omega,e^{-\varphi}dV)\longrightarrow A_\varphi^2(\Omega),
 \quad
 (P_{\Omega,\varphi}u)(z)
 =\int_\Omega K_{\Omega,\varphi}(z,w)
 u(w)e^{-\varphi(w)}\,dV(w).
\]
For the zero weight, we write
\[
 A^2(\Omega):=A_0^2(\Omega),\quad
 K_\Omega:=K_{\Omega,0},\quad
 P_\Omega:=P_{\Omega,0}.
\]
In particular, for $\varphi_t=-t\log\delta_\Omega$, $t>-1$,
we use the notation
\[
 A_t^2(\Omega):=A_{\varphi_t}^2(\Omega),\quad
 K_{\Omega,t}:=K_{\Omega,\varphi_t},\quad
 P_{\Omega,t}:=P_{\Omega,\varphi_t},
\]
and
\[
 \|f\|_t^2:=\|f\|_{\varphi_t}^2
 =\int_\Omega |f|^2\delta_\Omega^t\,dV.
\]

Assume that $\partial\Omega$ is $C^2$. For 
$\varepsilon>0$, let $\Omega_\varepsilon=\{\delta_\Omega>\varepsilon\}$.
The holomorphic Hardy space is
\[
 H^2(\Omega)=\left\{f\in\mathcal O(\Omega):
 \|f\|_{H^2(\Omega)}^2:= \limsup_{\varepsilon\downarrow0}
 \int_{\partial \Omega_\varepsilon}|f|^2\,d\sigma_\varepsilon<\infty\right\}.
\]
Its reproducing kernel $S_\Omega$ satisfies
\[
S_\Omega(z)=\sup_{0\ne f\in H^2(\Omega)}
 \frac{|f(z)|^2}{\|f\|_{H^2(\Omega)}^2}.
\]
By \cite{Stein}, every $f\in H^2(\Omega)$ has non-tangential
boundary values $f^*\in L^2(\partial\Omega)$, and
\begin{align*}
 \|f\|_{H^2(\Omega)}^2
 =\int_{\partial\Omega}|f^*|^2\,d\sigma,
\ \ \
 f(z)=\int_{\partial\Omega}S_\Omega(z,\zeta)f^*(\zeta)\,d\sigma(\zeta).
\end{align*}

We recall the following H\"ormander-type $L^2$ estimate.
\begin{theorem}[{\cite[Theorem~$1.1'$]{Chen14}}]\label{lem:chen-solution}
Let $\Omega\subset\C^n$ be a bounded pseudoconvex domain with $C^2$
boundary, and let $\hat\Omega\subset\Omega$ be a pseudoconvex subdomain.
Let
$\varphi\in\operatorname{PSH}(\hat\Omega)$ satisfy
$i\partial\bar\partial\varphi\ge i\partial\bar\partial\psi$ in the sense of
distributions, where $\psi\in C^2(\hat\Omega)\cap\operatorname{PSH}(\hat\Omega)$.
Then, for each $t>-1$ and each $\bar\partial$-closed $(0,1)$-form $v$ with
\[
 \int_{\hat\Omega}|v|^2_{i\partial\bar\partial\psi}
 e^{-\varphi}\delta_\Omega^t\,dV<\infty,
\]
there is a solution $\bar\partial u=v$ on $\hat\Omega$ such that
\begin{equation}\label{pre:ambient-solution}
 \int_{\hat\Omega}|u|^2e^{-\varphi}\delta_\Omega^t\,dV
 \le C_{t,\Omega}\int_{\hat\Omega}
 |v|^2_{i\partial\bar\partial\psi}e^{-\varphi}\delta_\Omega^t\,dV.
\end{equation}
\end{theorem}
The local uniformity of $C_{t,\Omega}$ for $t$ in compact subintervals
of $(-1,\infty)$ is implicit in Chen's proof
\cite[Section~2]{Chen14}; see also
\cite[Remark~2.5]{MX}.

\subsection{Green functions and barriers}\label{sec:pluripotential-prelim}
Assume that
$\Omega$ is bounded and hyperconvex, that is, it admits a negative
continuous plurisubharmonic exhaustion.
Every bounded pseudoconvex domain with $C^2$ boundary is
hyperconvex \cite{DF}.
Write $g_\Omega(z,w)$ for the negative
pluricomplex Green function and $g_{\Omega,w}:=g_\Omega(\cdot,w)$, namely,
\[
 g_\Omega(z,w)=\sup\{u(z):u\in\operatorname{PSH}(\Omega),\ u<0,
 \ u(\zeta)\le\log|\zeta-w|+O(1)\text{ near }w\}.
\]
For a fixed pole, it is continuous away from the pole, has zero
boundary values, and has relatively compact negative sublevel sets;
see \cite[Proposition~3.3.1 and Theorem~3.3.4]{BlockiNotes}.
For a fixed pole $w$, write
\[
 \Omega_g=\{z\in\Omega:g_\Omega(z,w)<-1\}
\]
for the Green sublevel set.
We also recall the Diederich--Forn\ae ss barriers.
\begin{theorem}[\cite{DF,CF}]
\label{lem:df-barriers}
Let $\Omega\subset\C^n$ be a bounded pseudoconvex domain with $C^2$ boundary.
For every $p\in\partial\Omega$ and $0<\eta<1$, there exist
a neighborhood $U$ of $p$, a $C^2$ pseudoconvex subdomain
$\hat\Omega\subset U\cap\Omega$ agreeing with $\Omega$ near $p$,
and a local $C^2$ defining function $\rho$ for $\Omega$ such that
\[
\psi=-(-\rho)^\eta\in\operatorname{PSH}(\hat\Omega),
\quad
-\psi\asymp\delta_\Omega^\eta
\quad\text{on }\hat\Omega.
\]
Moreover, $\delta_{\hat\Omega}=\delta_\Omega$ near $p$.
\end{theorem}
\section{Proof of the weighted comparison}\label{sec:two-sided}
From this section onward, $\Omega\Subset\C^n$, $n\ge2$, denotes a bounded
pseudoconvex domain with $C^2$ boundary. For this domain and its subdomains,
write
\[
L_\Omega(z):=|\log\delta_\Omega(z)|,\quad
L_{\hat\Omega}(z):=|\log\delta_{\hat\Omega}(z)|
\]
when the corresponding boundary distance is less than $e^{-1}$.
We first combine the Green-function cutoff argument in
\cite[Section~2]{Blocki25} with the distance-weighted estimate of
Theorem~\ref{lem:chen-solution}.
\begin{lemma}\label{lem:green}
Let $\hat\Omega\subset\Omega$ be a pseudoconvex subdomain with $C^2$
boundary, and let $w\in\hat\Omega$. For $t>-1$ and
$f\in\mathcal O(\hat\Omega_g)$ there is $F\in\mathcal O(\hat\Omega)$ such that
$F(w)=f(w)$ and
\begin{equation}\label{pre:extension}
 \int_{\hat\Omega}|F|^2\delta_\Omega^t\,dV
 \le C_{t,\Omega,n}\int_{\hat\Omega_g}|f|^2\delta_\Omega^t\,dV.
\end{equation}
The constants are uniform for $t$ in compact subintervals of $(-1,\infty)$.
\end{lemma}
\begin{proof}
If the right-hand side is infinite, take $F\equiv f(w)$. Otherwise,
fix $I\Subset(-1,\infty)$ and $t\in I$, and put
$
 g=g_{\hat\Omega,w}.
$
Choose a pseudoconvex exhaustion
$w\in\hat\Omega_1\Subset\hat\Omega_2\Subset\cdots\Subset\hat\Omega$.
By convolution regularization, choose
$g_\nu\in C^\infty(\hat\Omega_\nu)
\cap\operatorname{PSH}(\hat\Omega_\nu)$ such that
\[
g\le g_\nu<0,\quad
g_\nu\to g
\quad\text{locally uniformly on }
\hat\Omega\setminus\{w\}.
\]
Choose $\chi\in C^\infty(\R,[0,1])$ with
$\chi|_{(-\infty,-3]}=1$ and $\chi|_{[-2,\infty)}=0$. Set
\[
 \varphi_\nu=2ng_\nu-\log(1-g_\nu)\le 0,\ \ \
 v_\nu=\bar\partial(\chi(g_\nu)f)
       =\chi'(g_\nu)f\bar\partial g_\nu,
\]
where $\chi(g_\nu)f$ is extended by zero outside $\hat\Omega_g$.
Then
\[
 \bar\partial v_\nu=0,\quad
 \operatorname{supp}v_\nu\subset\{-3\le g_\nu\le-2\}
  \subset\hat\Omega_g,\quad
  -\varphi_\nu\le 6n+\log4\quad\text{on }\operatorname{supp}v_\nu,
\]
and
\begin{align*}
 i\partial\bar\partial\varphi_\nu
 =\left(2n+\frac1{1-g_\nu}\right)
       i\partial\bar\partial g_\nu
   +\frac{i\partial g_\nu\wedge\bar\partial g_\nu}{(1-g_\nu)^2}\ge \frac{i\partial g_\nu\wedge\bar\partial g_\nu}{(1-g_\nu)^2}.
\end{align*}
Theorem~\ref{lem:chen-solution} supplies a solution $u_\nu$ on
$\hat\Omega_\nu$ with $\bar\partial u_\nu=v_\nu$ and
\[
 \int_{\hat\Omega_\nu}|u_\nu|^2e^{-\varphi_\nu}\delta_\Omega^t\,dV
 \le C_{I,\Omega,n}\int_{\hat\Omega_\nu}|v_\nu|^2_{i\partial\bar\partial\varphi_\nu}e^{-\varphi_\nu} \delta_\Omega^t\,dV
 \le C_{I,\Omega,n}\int_{\hat\Omega_g}|f|^2\delta_\Omega^t\,dV.
\]
Set $F_\nu=\chi(g_\nu)f-u_\nu\in\mathcal O(\hat\Omega_\nu)$. Then
\begin{align*}
 \int_{\hat\Omega_\nu}|F_\nu|^2\delta_\Omega^t\,dV
 \le2\int_{\hat\Omega_g}|f|^2\delta_\Omega^t\,dV+2\int_{\hat\Omega_\nu}|u_\nu|^2
                   e^{-\varphi_\nu}\delta_\Omega^t\,dV
 \le C_{I,\Omega,n}\int_{\hat\Omega_g}|f|^2\delta_\Omega^t\,dV.
\end{align*}
After passing to a subsequence, which we do not relabel,
$F_\nu\to F\in\mathcal O(\hat\Omega)$ locally uniformly. Put
$\varphi=2ng-\log(1-g)$ and $u=\chi(g)f-F$. Fatou's lemma yields
\[
 \int_{\hat\Omega}|F|^2\delta_\Omega^t\,dV\le C_{I,\Omega,n}\int_{\hat\Omega_g}|f|^2\delta_\Omega^t\,dV,
 \quad
 \int_{\hat\Omega}|u|^2e^{-\varphi}\delta_\Omega^t\,dV
 \le C_{I,\Omega,n}\int_{\hat\Omega_g}|f|^2\delta_\Omega^t\,dV.
\]
For sufficiently small $\varepsilon>0$, on $B(w,\varepsilon)$ one has
\[
 \chi(g)=1,\quad u=f-F\in\mathcal O(B(w,\varepsilon)),\quad
 e^{-\varphi}=(1-g)e^{-2ng}\ge c_w|\zeta-w|^{-2n}.
\]
Consequently,
\[
 \int_{B(w,\varepsilon)}|u(\zeta)|^2|\zeta-w|^{-2n}\,dV(\zeta)<\infty
 \quad\Longrightarrow\quad u(w)=0,
\]
so $F(w)=f(w)$ and \eqref{pre:extension} follows.
\end{proof}

\begin{lemma}\label{lem:weighted-localization}
Let \(\hat\Omega\subset\Omega\) be a subdomain agreeing with \(\Omega\) near some \(p\in\partial\Omega\), and let $I\Subset(-1,\infty)$.
For $w\in\hat\Omega$ near $p$,
$t\in I$, and $f\in A_t^2(\hat\Omega)$, there exists
$F\in A_t^2(\Omega)$ such that
\[
 F(w)=f(w),\quad
 \int_\Omega|F|^2\delta_\Omega^t\,dV
 \le C_I\int_{\hat\Omega}|f|^2\delta_{\hat\Omega}^t\,dV.
\]
\end{lemma}
\begin{proof}
Choose neighborhoods $V\Subset U$ of $p$ such that
$\hat\Omega\cap U=\Omega\cap U$, and take
$\chi\in C_c^\infty(U)$ equal to $1$ near $\overline V$, with support
small enough that $\delta_{\hat\Omega}=\delta_\Omega$ on
$\hat\Omega\cap\operatorname{supp}\chi$.
Then
\[
c:=\operatorname{dist}
(\overline V,\operatorname{supp}\bar\partial\chi)>0,\quad
|\zeta-w|\ge c
\quad(w\in V\cap\hat\Omega,\
\zeta\in\operatorname{supp}\bar\partial\chi).
\]
Extend $\chi f$ by zero to $\Omega$; its distributional
$\bar\partial$ is $f\bar\partial\chi$, also extended by zero.
Apply Theorem~\ref{lem:chen-solution} with
\[
v=f\bar\partial\chi,\quad
\varphi(\zeta)=2n\log|\zeta-w|+A|\zeta|^2,
\quad \psi(\zeta)=A|\zeta|^2,
\]
where $A>0$ is fixed. Since
$i\partial\bar\partial\varphi\ge i\partial\bar\partial\psi$
and $|\zeta-w|\ge c$ on $\operatorname{supp}\bar\partial\chi$,
uniformly for $w\in V\cap\hat\Omega$, we obtain
$\bar\partial u=v$ with
\[
u(w)=0,\quad
\int_\Omega|u|^2\delta_\Omega^t\,dV
\le C_I
\int_{\hat\Omega\cap\operatorname{supp}\bar\partial\chi}
|f|^2\delta_\Omega^t\,dV.
\]
The logarithmic pole gives $u(w)=0$, as in Lemma~\ref{lem:green}.
Thus $F=\chi f-u$ is holomorphic on $\Omega$, has value
$f(w)$ at $w$, and satisfies
\[
 \int_\Omega|F|^2\delta_\Omega^t\,dV
 \le C_I\int_{\hat\Omega\cap\operatorname{supp}\chi}
              |f|^2\delta_\Omega^t\,dV
     =C_I\int_{\hat\Omega\cap\operatorname{supp}\chi}
              |f|^2\delta_{\hat\Omega}^t\,dV
              \le C_I\int_{\hat\Omega}|f|^2\delta_{\hat\Omega}^t\,dV.
\]
\end{proof}
The same argument applies with $\delta_\Omega^t$ in place of
$\delta_{\hat\Omega}^t$ on the right-hand side, provided that
$f\in\mathcal O(\hat\Omega)$ and
$
\int_{\hat\Omega}|f|^2\delta_\Omega^t\,dV<\infty.
$
Moreover, for $f\in A_t^2(\hat\Omega)$, the extension
$F\in A_t^2(\Omega)$ satisfies $F(w)=f(w)$ and
\[
|f(w)|^2=|F(w)|^2
\le K_{\Omega,t}(w)\|F\|_{A_t^2(\Omega)}^2
\le C_I K_{\Omega,t}(w)
       \|f\|_{A_t^2(\hat\Omega)}^2.
\]
Taking the supremum over $f\ne0$ gives the first inequality below.
The case $t=0$ gives the second. Thus, for $w$ near $p$,
\begin{equation}\label{local:weighted}
K_{\Omega,t}(w)\ge c_IK_{\hat\Omega,t}(w),
\quad
K_{\hat\Omega}(w)\le CK_\Omega(w).
\end{equation}

\subsection{The upper bound in Theorem~\ref{thm:comparison}}\label{sec:upper}
Let $W^s(\Omega)$ denote the $L^2$ Sobolev space of order
$s=0,1$, with $W^0(\Omega)=L^2(\Omega)$.
Let
\[
\tau:W^1(\Omega)\cap\mathcal O(\Omega)
\longrightarrow L^2(\partial\Omega),
\quad \tau f=f^*
\]
denote the Hardy boundary trace operator. Then
$
\|f\|_{H^2(\Omega)}
=\|\tau f\|_{L^2(\partial\Omega)}.
$
The multiplicative Sobolev trace inequality
\cite[Theorem~1.5.1.10]{Grisvard85} yields
\[
\|\tau f\|_{L^2(\partial\Omega)}^2
\le C\|f\|_{L^2(\Omega)}\|f\|_{W^1(\Omega)},
\quad f\in W^1(\Omega)\cap\mathcal O(\Omega).
\]
Kohn established weighted Sobolev regularity for smooth pseudoconvex
boundaries \cite{Kohn73}; Harrington
\cite[Corollary~1(4)]{Harrington09} supplies the $C^2$ case used here.
Fix $A>0$ sufficiently large and put $\varphi(\zeta)=A|\zeta|^2$.
The projection $P_{\Omega,\varphi}$ introduced in
Section~\ref{sec:holomorphic-prelim} is bounded on $W^1(\Omega)$.
Equivalence of the weighted and unweighted $L^2$ norms also gives its
boundedness on $L^2(\Omega)$.
Combining these bounds, we obtain, for $u\in W^1(\Omega)$,
\begin{equation}\label{eq:projection-trace}
\begin{aligned}
\|P_{\Omega,\varphi} u\|_{H^2(\Omega)}^2
=\|\tau (P_{\Omega,\varphi} u)\|_{L^2(\partial\Omega)}^2
\le C\|P_{\Omega,\varphi} u\|_{L^2(\Omega)}
       \|P_{\Omega,\varphi} u\|_{W^1(\Omega)}
\le C\|u\|_{L^2(\Omega)}\|u\|_{W^1(\Omega)}.
\end{aligned}
\end{equation}
\begin{lemma}
\label{lem:weighted-decomposition}
For every $-1<t_0<0$, there exist $r_0,C>0$
such that, for every $t\in[t_0,0)$,
$f\in A_t^2(\Omega)$, and $0<r<r_0$, there exists
$u_r\in A^2(\Omega)$ satisfying
$f-u_r\in H^2(\Omega)$ and
\begin{equation}\label{direct:decomposition}
\|u_r\|_{A^2(\Omega)}
\le Cr^{-t/2}\|f\|_t,\quad
\|f-u_r\|_{H^2(\Omega)}
\le Cr^{-(1+t)/2}\|f\|_t.
\end{equation}
\end{lemma}

\begin{proof}
Fix $-1<t_0<0$, $t\in[t_0,0)$, and $f\in A_t^2(\Omega)$. Since
$A_t^2(\Omega)\subset A^2(\Omega)$, we have
$P_{\Omega,\varphi}f=f$.
Choose $0<r_0<1$ sufficiently small that $\delta_\Omega$ is $C^2$
on $\{0<\delta_\Omega<2r_0\}$, and take $\eta\in C^\infty([0,\infty))$ such that
\[
0\le\eta\le1,\quad
\eta=1\ \text{on }[0,1],\quad
\eta=0\ \text{on }[2,\infty).
\]
For $0<r<r_0$, set
\[
\chi_r=\eta(\delta_\Omega/r),\quad
u_r=P_{\Omega,\varphi}(\chi_rf),\quad
v_r=f-u_r=P_{\Omega,\varphi}((1-\chi_r)f).
\]
Then $u_r\in A^2(\Omega)$ and
\[
\begin{aligned}
\|u_r\|_{A^2(\Omega)}^2
\le C\|\chi_rf\|_{L^2(\Omega)}^2
\le C\int_{\{\delta_\Omega<2r\}}|f|^2\,dV
\le Cr^{-t}\|f\|_t^2.
\end{aligned}
\]
Choose $j_0\ge0$ such that
$r_{j_0}=2^{j_0}r\in(r_0/2,r_0]$, and put
\[
r_j=2^jr,\quad
\psi_j=\chi_{r_{j+1}}-\chi_{r_j}
\quad(0\le j<{j_0}),\quad
\psi_{j_0}=1-\chi_{r_{j_0}}.
\]
Thus
\[
1-\chi_r=\sum_{j=0}^{j_0}\psi_j,\quad
v_r=\sum_{j=0}^{j_0}P_{\Omega,\varphi}(\psi_jf).
\]
Let
\(
U_j=\{r_j/2<\delta_\Omega<8r_j\}
\ (j<{j_0}),\
U_{j_0}=\{\delta_\Omega>r_{j_0}/2\},
\)
and set
\[
c_j^2=\int_{U_j}|f|^2\delta_\Omega^t\,dV.
\]
Since
\(
{\delta_\Omega(z)}/8<r_j<2\delta_\Omega(z)
\) for \(z\in U_j,\ j<{j_0}\),
we have
\[
\sum_{j=0}^{j_0}\mathbf1_{U_j}\le5.
\]
Integration therefore gives
\[
\begin{aligned}
\sum_{j=0}^{j_0}c_j^2
=\int_\Omega |f|^2\delta_\Omega^t
\sum_{j=0}^{j_0}\mathbf1_{U_j}\,dV
\le5\|f\|_t^2.
\end{aligned}
\]
By the construction of $\psi_j$ and $U_j$,
\[
\operatorname{supp}\psi_j\subset U_j,\quad
|\psi_j|\le1,\quad
|\nabla\psi_j|\le Cr_j^{-1}.
\]
Since
$B(z,r_j/4)\subset U_j$ for
$z\in\operatorname{supp}\psi_j$,
the $L^2$ Cauchy estimates and Fubini's theorem yield
\[
\begin{aligned}
\int_{\operatorname{supp}\psi_j}|\nabla f|^2\,dV
\le Cr_j^{-2n-2}
\int_{\operatorname{supp}\psi_j}
\int_{B(z,r_j/4)}|f(\zeta)|^2\,dV(\zeta)\,dV(z)
\le Cr_j^{-2}\int_{U_j}|f|^2\,dV.
\end{aligned}
\]
Since $\delta_\Omega\asymp r_j$ on $U_j$ for $j<j_0$,
while $r_{j_0}\asymp r_0$ and
$U_{j_0}=\{\delta_\Omega>r_{j_0}/2\}$, we have
\[
\|\psi_jf\|_{L^2(\Omega)}^2
\le
\int_{U_j}|f|^2\,dV
\le Cr_j^{-t}c_j^2,\quad 0\le j\le {j_0}.
\]
The product rule now gives
\[
\begin{aligned}
\|\psi_jf\|_{W^1(\Omega)}^2
\le C\int_{U_j}|f|^2\,dV
 +C\int_{\operatorname{supp}\psi_j}
       \bigl(|\nabla f|^2+r_j^{-2}|f|^2\bigr)\,dV
\le Cr_j^{-2-t}c_j^2.
\end{aligned}
\]
Hence \eqref{eq:projection-trace} gives
\[
\begin{aligned}
\|P_{\Omega,\varphi}(\psi_jf)\|_{H^2(\Omega)}^2
\le C\|\psi_jf\|_{L^2(\Omega)}
       \|\psi_jf\|_{W^1(\Omega)}
\le Cr_j^{-(1+t)}c_j^2.
\end{aligned}
\]
Summing in $j$, we find that $v_r\in H^2(\Omega)$ and
\[
\begin{aligned}
\|v_r\|_{H^2(\Omega)}
\le C\sum_{j=0}^{j_0}r_j^{-\frac{1+t}{2}}c_j
\le C\left(\sum_{j=0}^{j_0}r_j^{-(1+t)}\right)^{1/2}
        \left(\sum_{j=0}^{j_0}c_j^2\right)^{1/2}
\le \frac{C}{\sqrt{1-2^{-(1+t_0)}}} r^{-\frac{1+t}{2}}\|f\|_t.
\end{aligned}
\]
This proves \eqref{direct:decomposition} uniformly for
$t\in[t_0,0)$.
\end{proof}
\begin{proposition}\label{prop:kernel-inequality}
For $-1<t<0$ and $z\in\Omega$,
\begin{equation}\label{fixed:kernel}
 K_{\Omega,t}(z)\le C_{\Omega,t}
 K_\Omega(z)^{1+t}S_\Omega(z)^{-t}.
\end{equation}
For every fixed $-1<t_0<0$, the constants can be chosen uniformly for
$t\in[t_0,0)$.
\end{proposition}
\begin{proof}
Fix $-1<t_0<0$, $t\in[t_0,0)$, and
$f\in A_t^2(\Omega)$.
By Lemma~\ref{lem:weighted-decomposition},
for every $0<r<r_0$, there exists $u_r\in A^2(\Omega)$
such that $f-u_r\in H^2(\Omega)$ and
\[
\begin{aligned}
|f(z)|
&\le |u_r(z)|+|(f-u_r)(z)|\\
&\le K_\Omega(z)^{1/2}\|u_r\|_{A^2(\Omega)}
   +S_\Omega(z)^{1/2}\|f-u_r\|_{H^2(\Omega)}\\
&\le C\left(
K_\Omega(z)^{1/2}r^{-t/2}
+S_\Omega(z)^{1/2}r^{-(1+t)/2}
\right)\|f\|_t.
\end{aligned}
\]
If $S_\Omega(z)/K_\Omega(z)<r_0$, taking
\(
r=S_\Omega(z)/K_\Omega(z)
\)
gives
\[
|f(z)|^2
\le C K_\Omega(z)^{1+t}S_\Omega(z)^{-t}\|f\|_t^2.
\]
Otherwise, putting $R=\sup_\Omega\delta_\Omega$, we have
\[
\begin{aligned}
|f(z)|^2
\le K_\Omega(z)\|f\|_{A^2(\Omega)}^2
\le R^{-t}K_\Omega(z)\|f\|_t^2
\le (R/r_0)^{-t}
K_\Omega(z)^{1+t}S_\Omega(z)^{-t}\|f\|_t^2.
\end{aligned}
\]
Taking the supremum over $f\ne0$ proves
\eqref{fixed:kernel}, uniformly for $t\in[t_0,0)$.
If $C_d$ denotes the constant in \eqref{direct:decomposition},
one may take the uniform coefficient
\[
 C_{\Omega,t_0}=\max\{4C_d^2,1,R/r_0\}.
\]
\end{proof}

\begin{proof}[Proof of the upper bound in
Theorem~\ref{thm:comparison}]
Fix $0<\alpha<1$ and $I\Subset(-1,\infty)$.
For $p\in\partial\Omega$, let $\hat\Omega$ and $\psi$
be as in Theorem~\ref{lem:df-barriers}, with
$\eta=\alpha$. For $w\to p$ in $\hat\Omega$, put
$r=\delta_{\hat\Omega}(w)=\delta_\Omega(w)$.

The Green-function comparison in \cite[(4.4)--(4.7)]{CF},
applied to $g_{\hat\Omega}$ and $\psi$, gives
\[
g_{\hat\Omega}(z,w)
\ge -C\left(\frac{\delta_\Omega(z)}r\right)^\alpha
L_\Omega(w),
\quad z\in\hat\Omega\setminus B(w,r/2).
\]
Since $\delta_\Omega(z)\ge r/2$ on $B(w,r/2)$,
we obtain
\begin{equation}\label{upper:green-inclusion}
\hat\Omega_g\subset
\left\{z\in\hat\Omega:
\delta_\Omega(z)\ge
crL_\Omega(w)^{-1/\alpha}\right\}.
\end{equation}
Applying Lemma~\ref{lem:green} to
\[
f(z)=\frac{K_{\hat\Omega}(z,w)}
{\sqrt{K_{\hat\Omega}(w)}},
\quad \|f\|_{A^2(\hat\Omega)}=1,
\]
and using \eqref{upper:green-inclusion}, we obtain
$h\in\mathcal O(\hat\Omega)$ satisfying
\[
h(w)=K_{\hat\Omega}(w)^{1/2},\quad
\int_{\hat\Omega}|h|^2\delta_\Omega^{-1/2}\,dV
\le Cr^{-1/2}L_\Omega(w)^{1/(2\alpha)}.
\]
The extension estimate following Lemma~\ref{lem:weighted-localization},
with $\delta_\Omega^{-1/2}$ on both sides, provides
$F\in A_{-1/2}^2(\Omega)$ such that $F(w)=h(w)=K_{\hat\Omega}(w)^{1/2}$ and
\[
\begin{aligned}
\int_\Omega |F|^2\delta_\Omega^{-1/2}\,dV
\le C\int_{\hat\Omega}|h|^2\delta_\Omega^{-1/2}\,dV
\le Cr^{-1/2}L_\Omega(w)^{1/(2\alpha)}.
\end{aligned}
\]
Therefore
\[
 K_{\Omega,-1/2}(w)
 \ge c r^{1/2}L_\Omega(w)^{-1/(2\alpha)}K_{\hat\Omega}(w)
 \ge c r^{1/2}L_\Omega(w)^{-1/(2\alpha)}K_\Omega(w),
 \begingroup
\renewcommand{\thefootnote}{\fnsymbol{footnote}}
\footnote[1]{This inequality slightly improves an estimate obtained by
Chen--Zhang in the proof of
\cite[Theorem~7.8, p.~55]{ChenZhang22}.}
\endgroup
\]
where the last inequality uses domain monotonicity:
$K_{\hat\Omega}(w)\ge K_\Omega(w)$.
On the other hand, Proposition~\ref{prop:kernel-inequality}, with
$t=-1/2$, gives
\[
 K_{\Omega,-1/2}(w)
 \le C K_\Omega(w)^{1/2}S_\Omega(w)^{1/2}.
\]
Combining these inequalities yields
\begin{equation}\label{fixed:upper}
 \frac{K_\Omega(w)}{S_\Omega(w)}
 \le C\frac{L_\Omega(w)^{1/\alpha}}r.
\end{equation}
For $-1<t<0$, Proposition~\ref{prop:kernel-inequality} yields
\[
 \frac{K_{\Omega,t}(w)}{S_\Omega(w)}\le C_I\left(\frac{K_\Omega(w)}{S_\Omega(w)}\right)^{1+t}
 \le C_I\left(\frac{L_\Omega(w)^{1/\alpha}}r\right)^{1+t}.
\]
The constants are uniform for $t\in I\cap(-1,0)$ by
Proposition~\ref{prop:kernel-inequality}; the case $t=0$ follows from
\eqref{fixed:upper}.
For $t\ge0$, Lemma~\ref{lem:green} with $t=0$,
\eqref{upper:green-inclusion}, and \eqref{local:weighted} give
\begin{align*}
 |f(w)|^2
 \le C K_{\hat\Omega}(w)
       \int_{\hat\Omega_g}|f|^2\,dV
 \le C_I K_\Omega(w)r^{-t}L_\Omega(w)^{t/\alpha}
                        \|f\|_{A_t^2(\Omega)}^2,
 \ \ \ \forall\,f\in A_t^2(\Omega).
\end{align*}
Consequently,
\[
 \frac{K_{\Omega,t}(w)}{S_\Omega(w)}
 \le C_Ir^{-t}L_\Omega(w)^{t/\alpha}\frac{K_\Omega(w)}{S_\Omega(w)}
 \le C_I\left(\frac{L_\Omega(w)^{1/\alpha}}r\right)^{1+t}.
\]
A finite boundary cover completes the proof.
\end{proof}

Taking $t=0$ in Theorem~\ref{thm:comparison} proves
Theorem~\ref{thm:unweighted}.

\section{Logarithmic obstructions}\label{sec:counterexamples}

The obstruction can be realized on a rotationally invariant Hartogs
domain in $\C^2$. Only its kernels on the normal axis are needed, so
the construction reduces to one-dimensional radial moments.

\subsection{The radial Hartogs domain and its kernels}
Let $q\in C^\infty([0,1])$ be nonnegative, nondecreasing, and flat
at zero, and set
\begin{equation}\label{examples:radial-profile}
 \phi(r)=\int_0^r\frac{q(u)}u\,du.
\end{equation}
Choose a smooth convex nondecreasing function
$\psi:[0,1)\to[0,\infty)$, vanishing on $[0,1/4]$, with
$\psi(u)=C_*-\frac12\log(1-u)$ near $1$. Put
\begin{equation}\label{examples:profile}
 \Phi(\zeta)=\phi(|\zeta|)+\psi(|\zeta|^2),\qquad
 \Omega_\Phi=\{(\zeta,w)\in\D\times\C:
 |w|<e^{-\Phi(\zeta)}\}.
\end{equation}
Then $\Omega_\Phi$ is a bounded $C^\infty$ pseudoconvex complete
Reinhardt domain, and $p=(0,1)$ is an infinite-type boundary point.
Indeed, $\Delta\phi(r)=q'(r)/r\ge0$ for $r>0$, while near
$|\zeta|=1$ the defining inequality becomes
$e^{2C_*+2\phi(|\zeta|)}|w|^2+|\zeta|^2<1$.

For $\lambda>0$, define
\[
 M(\lambda)=\int_\D e^{-2\lambda\Phi}\,dV,\qquad
 \widetilde M(\lambda)=\int_\D e^{-2\lambda\Phi}
 \sqrt{1+e^{-2\Phi}|\nabla\Phi|^2}\,dV.
\]
Radial dilation on this strictly star-shaped domain approximates $H^2$
functions by functions holomorphic
on a neighborhood of $\overline\Omega_\Phi$; their Taylor expansions
show that the monomials form a complete orthogonal system.
Rotational orthogonality and direct integration give
\[
 \|w^\nu\|_{A^2(\Omega_\Phi)}^2
 =\frac{\pi}{\nu+1}M(\nu+1),\qquad
 \|w^\nu\|_{H^2(\Omega_\Phi)}^2
 =2\pi\widetilde M(\nu+1/2).
\]
All monomials containing a positive power of $\zeta$ vanish on the
normal axis. Thus the Forelli--Rudin expansions reduce to
\begin{equation}\label{examples:kernel-series}
 \begin{aligned}
 K_{\Omega_\Phi}(0,w)
 &=\frac1\pi\sum_{\nu\ge0}
   \frac{(\nu+1)|w|^{2\nu}}{M(\nu+1)},\\
 S_{\Omega_\Phi}(0,w)
 &=\frac1{2\pi}\sum_{\nu\ge0}
   \frac{|w|^{2\nu}}{\widetilde M(\nu+1/2)};
 \end{aligned}
\end{equation}
see also \cite[(33)--(34)]{EnglisXu15}.

\subsection{The plateau construction}
\begin{proposition}\label{ex:obstructions-1}
Let $H:[1,\infty)\to(0,\infty)$ be nondecreasing and satisfy
$H(s)=o(s)$. There is a function $\Phi$ of the form
\eqref{examples:profile} such that $\Omega=\Omega_\Phi$ admits
a sequence $z_k\to p=(0,1)$ along the inward normal with
\begin{equation}\label{examples:upper-limit}
 \frac{\delta_\Omega(z_k)K_\Omega(z_k)}
 {S_\Omega(z_k)H(L_\Omega(z_k))}\longrightarrow\infty.
\end{equation}
\end{proposition}

\begin{proof}
Choose integers $T_k$ recursively so that $T_1=10$ and
\begin{equation}\label{examples:scale-growth}
 T_{k+1}\ge T_k+1+\max\{T_k^2,e^{k^2}\},\qquad
 H(2kT_{k+1})\le T_{k+1}/k^2.
\end{equation}
Such a choice is possible because $H(s)/s\to0$. Set
\begin{equation}\label{examples:scales}
 G_k=e^{-kT_{k+1}},\qquad N_k=T_{k+1}-T_k-1.
\end{equation}
Fix a nondecreasing $\chi\in C^\infty(\R)$ with
$\chi=0$ on $(-\infty,0]$ and $\chi=1$ on $[1,\infty)$, and put
\begin{equation}\label{upper:profile-choice}
 q(0)=0,\qquad
 q(r)=\sum_{\ell\ge1}(G_\ell-G_{\ell+1})
 \chi(\log r+T_{\ell+1}),\quad 0<r\le1.
\end{equation}
The coefficients are positive. Moreover, for fixed
$m,M\in\mathbb N_0$,
\[
 \sup_{e^{-T_{k+1}}\le r\le e^{-T_k}}
 r^{-M}|q^{(m)}(r)|
 \le C_{m,M}e^{-(k-m-M)T_{k+1}}\longrightarrow0.
\]
Hence $q$ is smooth, nondecreasing, and flat at zero.

Fix a sufficiently large $k$ and abbreviate
\[
 G=G_k,\quad N=N_k,\quad b=e^{-T_k},\quad a=be^{-N}.
\]
The definition of $q$ gives
\begin{equation}\label{examples:plateau-data}
 \begin{gathered}
 q=G\ \hbox{on }[a,b],\qquad
 q\ge e^{-(k-1)T_k}\ge4G\ \hbox{on }[eb,1],\\
 \phi(a)\le G+\sum_{\ell>k}T_{\ell+1}e^{-\ell T_{\ell+1}}
 \le CG,\qquad
 \phi(b)=\phi(a)+NG\asymp NG.
 \end{gathered}
\end{equation}
Let $w_k=e^{-\phi(b)}$ and $z_k=(0,w_k)$. Since $N\asymp T_{k+1}$
and $NG\to0$, the tubular neighborhood theorem gives
\begin{equation}\label{examples:normal-depth}
 \delta_\Omega(z_k)=1-w_k\asymp NG,\qquad
 L_\Omega(z_k)\le2kT_{k+1}
\end{equation}
for all sufficiently large $k$.

For $x=\lambda G\in[1/2,3/4]$, the identity
$\phi(r)=\phi(b)+G\log(r/b)$ on $[a,b]$ gives
\begin{equation}\label{examples:annular-moment}
 \int_a^b r e^{-2\lambda\phi(r)}\,dr
 =b^2e^{-2\lambda\phi(b)}
 \frac{1-e^{-2N(1-x)}}{2(1-x)}.
\end{equation}
The remaining pieces obey
\[
 \begin{aligned}
 \frac{\int_0^a r e^{-2\lambda\phi(r)}\,dr}
 {b^2e^{-2\lambda\phi(b)}}
 &\le\tfrac12e^{-2N(1-x)+2\lambda\phi(a)}\le C,\\
 \frac{\int_b^1 r e^{-2\lambda\phi(r)}\,dr}
 {b^2e^{-2\lambda\phi(b)}}
 &\le\int_0^1e^{2(1-x)u}\,du
 +e^{2(1-x)}\int_0^\infty e^{2(1-4x)u}\,du\le C.
 \end{aligned}
\]
Since $\psi\ge0$, this proves
\begin{equation}\label{examples:full-moments}
 M(\lambda)\le Cb^2e^{-2\lambda\phi(b)},
 \qquad \tfrac12\le\lambda G\le\tfrac34.
\end{equation}
The set of integers $\nu\ge0$ with
$\tfrac12\le(\nu+1)G\le\tfrac34$ has cardinality comparable to
$G^{-1}$, and each satisfies $\nu+1\asymp G^{-1}$. Therefore
\eqref{examples:kernel-series} and \eqref{examples:full-moments} yield
\[
 K_\Omega(z_k)\ge\frac{c}{b^2}
 \sum_{\frac12\le(\nu+1)G\le\frac34}(\nu+1)
 \ge\frac{c}{G^2b^2}.
\]

For the Szeg\H{o} kernel, put $y_\nu=(\nu+1/2)G$.
Since $\psi=0$ on $[a^2,b^2]$, the same annular integral gives
\[
 \widetilde M(\nu+1/2)\ge M(\nu+1/2)
 \ge\pi b^2w_k^{2\nu+1}
 \frac{1-e^{-2N(1-y_\nu)}}{1-y_\nu},
\]
where the last quotient has value $2N$ at $y_\nu=1$.
For $N\ge1$,
\[
 \frac{1-y}{1-e^{-2N(1-y)}}\le
 \begin{cases}
 C,&0\le y\le1,\\
 C N^{-1}e^{-N(y-1)},&y>1.
 \end{cases}
\]
Consequently, as $w_k\to1$ and $NG\to0$,
\[
 \begin{aligned}
 S_\Omega(z_k)
 &\le\frac{C}{b^2}\sum_{\nu\ge0}
 \frac{1-y_\nu}{1-e^{-2N(1-y_\nu)}}\\
 &\le\frac{C}{b^2}\left(G^{-1}
 +\frac1{N(1-e^{-NG})}\right)
 \le\frac{C}{Gb^2}.
 \end{aligned}
\]
Combining these estimates with \eqref{examples:scale-growth} and
\eqref{examples:normal-depth}, we obtain
\[
 \frac{\delta_\Omega(z_k)K_\Omega(z_k)}
 {S_\Omega(z_k)H(L_\Omega(z_k))}
 \ge\frac{cN_k}{H(2kT_{k+1})}
 \ge ck^2\longrightarrow\infty.
\]
\end{proof}

\begin{remark}\label{rem:hartogs-powers}
The functions below may be modified on a bounded interval to be positive
and nondecreasing. Proposition~\ref{ex:obstructions-1} rules out general
upper estimates
\[
 \frac{K_\Omega(z)}{S_\Omega(z)}
 \le\frac{C_\Omega H(L_\Omega(z))}{\delta_\Omega(z)}
\]
on smoothly bounded pseudoconvex domains in $\C^2$, in particular for
\[
 H(s)=s^\vartheta\quad(0\le\vartheta<1),\qquad
 H(s)=\frac{s}{[\log(e+s)]^\gamma}\quad(\gamma>0).
\]
For these two families, a single domain suffices and the scales can
be chosen explicitly:
\[
 T_1=10,\qquad T_{k+1}=\lceil e^{T_k}\rceil,\qquad
 G_k=e^{-kT_{k+1}}.
\]
Indeed, the proof gives
$\delta_\Omega(z_k)K_\Omega(z_k)/S_\Omega(z_k)\ge cT_{k+1}$ and
$L_\Omega(z_k)\le2kT_{k+1}$. Hence the corresponding quotients
in \eqref{examples:upper-limit} are bounded below, respectively, by
\[
 c\frac{T_{k+1}^{1-\vartheta}}{k^\vartheta}
 \longrightarrow\infty,\qquad
 c\frac{[\log(e+2kT_{k+1})]^\gamma}{k}
 \longrightarrow\infty.
\]
For an arbitrary prescribed $H=o(s)$, the recursive choice in
\eqref{examples:scale-growth} gives the full statement above.
\end{remark}

\section*{Declaration of AI assistance}
The author used ChatGPT (OpenAI) to assist with the development and checking
of proofs, particularly the counterexample construction in
Section~\ref{sec:counterexamples}, and with writing, revision, and proofreading. The author reviewed and edited all AI-generated content
and takes full responsibility
for the content of this article.
\raggedbottom

\end{document}